\magnification=1200
\def\hal{{\vrule height 10pt width 4pt depth 0pt}}
\def\A{{\cal A}}
\def\C{{\bf C}}
\def\M{{\cal M}}
\def\N{{\bf N}}
\def\R{{\bf R}}

\centerline{\bf Locally nonconical convexity}
\medskip

\centerline{C.\ A.\ Akemann, G.\ C.\ Shell, and N.\ Weaver
\footnote{}{\noindent Keywords: convex set, locally nonconical, continuous section, continuous
selection, strictly convex, zonoid\break
2000 Mathematics Subject Classification 52A20, 52A07, 46A55, 46C05}}
\bigskip
\bigskip

{\it There is a hierarchy of structure conditions for convex sets. In this
paper we study a recently defined [3, 8, 9] condition called locally
nonconical convexity (abbreviated LNC).  Is is easy to show that every
strictly convex set is LNC, as are half-spaces and finite intersections of
sets of either of these types, but many more sets are LNC.  For instance,
every zonoid (the range of a nonatomic vector-valued measure) is LNC
(Corollary 34).  However, there are no infinite-dimensional compact LNC
sets (Theorem 23).
\medskip

The LNC concept originated in a search for continuous sections, and the
present paper shows how it leads naturally (and constructively) to
continuous sections in a variety of situations. Let $Q$ be a compact,
convex set in $\R^n$, and let $T$ be a linear map from $\R^n$ into
$\R^m$.  We show (Theorem 1) that $Q$ is LNC if and only if the
restriction of any such $T$ to $Q$ is an open map of $Q$ onto $T(Q)$.
This implies that if $Q$ is LNC, then any such $T$ has continuous sections
(i.e.\ there are continuous right inverses of $T$) that map from $T(Q)$ to
$Q$, and in fact it is possible to define continuous sections constructively
in various natural ways (Theorem 3, Corollary 4, and Theorem 5). If $Q$ is
strictly convex and $T$ is not 1-1, we can construct continuous sections
which take values in the boundary of $Q$ (Theorem 6).
\medskip

When we give up compactness it is natural to consider a closed, convex,
LNC subset $Q$ of a Hilbert space $X$ which may be infinite-dimensional.
In this case we must assume that $T$ is left Fredholm, i.e.\ a bounded linear
map with closed range and finite-dimensional kernel. We can then prove
results analogous to those mentioned in the last paragraph (Theorems
16-20). We also prove that $T(Q)$ is LNC (Theorem 25).}
\bigskip
\bigskip

\noindent {\bf 1. Introduction} \bigskip
The concept of a locally nonconical (abbr.\ LNC) convex set arose
in [3] and was explored further in the Ph.D. dissertation [8] of
the second author; some of the results in [8] appeared in [9].
The LNC concept originated in a search for continuous sections,
and the present paper shows why it is indeed the key to the
existence of continuous sections in a variety of situations.
\bigskip

\noindent {\bf Definition.} A convex set  $Q$ in a Hausdorff
topological vector space is {\it locally nonconical} (LNC) if for
any net $(x_t)$ in $Q$ that converges to a point $x$ in $Q$ and
for any other point $x'$ in $Q$, the points $x_t + (x' - x)/2$
eventually lie in $Q$.
\medskip

Equivalently, for every pair of points $x$ and $x'$
in $Q$ there exists a relative neighborhood  $U$ of $x$ in $Q$ such
that $U + (x'-x)/2 \subset Q$. Also, observe that if $Q$ is
LNC, then because $x_t + (x' - x)/2 \to (x' + x)/2$, we can apply
the LNC condition again with $(x' + x)/2$ in place of $x$, and conclude
that $x_t + {3\over 4}(x' - x)$ is eventually in $Q$, and so on;
inductively, the LNC condition implies that
for any $x, x', x_t \in Q$ with $x_t \to x$, we have
$x_t + e(x' - x) \in Q$ eventually, for any $e \in (0, 1)$.
\medskip

The name arose from the fact that no
point on the ``slanted'' portion of an ordinary circular cone can be
the point $x$ of this definition. See Example 9.
\medskip

In practice, to verify the LNC property it is sufficient to show that
there is a subnet of $x_t + (x' - x)/2$ which eventually lies in $Q$.
That is because if $x, x', x_t$ falsifies LNC, then by passing to a
subnet we can ensure that $x_t + (x' - x)/2$ is never in $Q$. We will
use this observation repeatedly.
\medskip

From here on we take, as part of the definition of topological vector
spaces, that they are Hausdorff. We also take the scalar field to be $\R$.
\medskip

There is a hierarchy of conditions that have
been defined for convex sets. It is shown in [8] that LNC lies below
``strictly convex'' and above ``uniformly stable'' in this hierarchy.
Indeed, if the net $(x_t)$ is not required to lie in $Q$, or if the
word ``relative'' is deleted from the second LNC definition, then the
condition becomes equivalent to strict convexity. What
sets do we get by slightly weakening the definition of ``strictly convex''?
Besides strictly convex sets, we get half spaces, finite intersections of
LNC sets, finite Cartesian products of LNC sets and all zonoids. See Sections
3 and 6.
\medskip

Since a strictly convex set must have nonvoid topological interior, it is
clear that there cannot be any compact, infinite-dimensional strictly convex
sets. The same holds for compact, infinite-dimensional LNC sets (see
Theorem 23).  The proof is more difficult in the LNC case because this
theorem, despite its negative sounding tone, is actually an existence
result! That is, we show that any infinite-dimensional, compact, convex
set $Q$ contains points $x$ and $x'$ and a net $\{x_t\}$  that converges
to $x$ such that  $x_t + (x' - x)/2$ never lies in $Q$. This implies that
there is a continous affine map out of $Q$ that is not an open map (see
Corollary 24). On the other hand, there are plenty of closed
infinite-dimensional LNC sets, and some of what we know about compact LNC
sets actually holds (with more difficult proofs) for closed
LNC sets; see Section 4.
\medskip

Here is a little more detail on our motivation for studying the LNC property.
Let $T: X \to Y$ be a linear map betwen locally convex topological vector
spaces. If $Q \subset X$ is convex, we regard $f = T|_Q: Q \to T(Q)$ as a
covering map and ask whether it possesses a continuous section. That is,
does there exist a continuous map $g: T(Q) \to Q$ such that $f \circ g =
{\rm id}_{T(Q)}$? Another way to put this is to ask if $f$ has a
continuous right inverse. Such a map $g$ is also sometimes called a
continuous selection, because for each $y \in T(Q)$ it selects,
in a continuous fashion, a point $g(y) \in f^{-1}(y)$.  Continuous
selections have been studied extensively (see e.g.\ [6]), but our approach
is different. We want a condition on the set $Q$ that will {\bf guarantee}
the existence of continuous sections for {\bf all} maps $f$ as described
above. We show in Section 2 that for compact sets $Q$ in $\R^n$ the LNC
condition fills the bill, and in Section 4 we treat the case of closed
sets.
\medskip

The LNC condition is not equivalent to the existence of continuous sections
(see Example 9). However, it seems that without LNC any continuous section
which does exist is somehow ``ad hoc,'' though it is difficult to make this
assertion precise. The principle we have in mind is that any reasonable general
construction of sections will automatically produce continuous sections if
and only if $Q$ is LNC. Again, what constitutes a ``reasonable'' construction is
unclear, but we will give several examples in Section 2.
\bigskip
\bigskip

\noindent {\bf 2. Continuous sections in finite dimensions}
\bigskip

Since our results are simplest and most complete in the finite-dimensional
case, we begin there. Let $T: \R^n \to \R^m$ be a linear map. If
$Q \subset \R^n$ is convex and compact, we define $f = T|_Q: Q \to T(Q)$
and ask when there exists a continuous map $g: T(Q) \to Q$ such that
$f \circ g = {\rm id}_{T(Q)}$.
\medskip

It follows from Michael's selection theorem [6, Examples 1.1, 1.1*, and
Theorem 3.2] that a continuous section exists if $f: Q \to T(Q)$ is an
open map. The LNC condition exactly guarantees that any restriction of a
linear map is
open. In addition, our methods prove more. Not only do continuous sections
always exist when $Q$ is LNC, but they can be defined constructively
and are of a minimizing type.
\bigskip

\noindent {\bf Theorem 1.} {\it Let $Q$ be a compact, convex subset of
$\R^n$. The following conditions are equivalent:
\medskip

{\narrower{
\noindent (i) $Q$ is LNC;
\medskip

\noindent (ii) for every linear map $T: \R^n \to \R^{n-1}$, the restriction
map $f = T|_Q: Q \to T(Q)$ is open;
\medskip

\noindent (iii) for every $m$ and every linear map $T: \R^n \to \R^m$,
the restriction map $f = T|_Q: Q \to T(Q)$ is open.
\medskip}}}

\noindent {\it Proof.} (i) $\Rightarrow$ (iii). Suppose $Q$ is LNC but
$f = T|_Q$ is not open for some $T: \R^n \to \R^m$. We shall reach a
contradiction. There exists $x' \in Q$ and a relative open neighborhood $U$
of $x'$ in $Q$ such that $f(U)$ is not a relative open neighborhood
of $f(x')$ in $T(Q)$. Thus, there is a sequence $(y_i)$ in $T(Q)$ which
converges to $f(x')$, such that $y_i$ is not in $f(U)$ for any $i$.
Find $(x_i) \subset Q$ such that $f(x_i) = y_i$
and pass to a subsequence so that $(x_i)$ converges.
Say $x_i \to x$; then
$$f(x) = \lim f(x_i) = \lim y_i = f(x').$$
Choose $e \in (0,1)$ small enough that $x' + e(x - x') \in U$; we can do this
because $U$ is relatively open and $x' + e(x - x') \in U$ for all $e \in [0,1]$.
Since $Q$ is LNC, $x_i + (1 - e)(x' - x)$ is
eventually in $Q$. But since $x_i + (1 - e)(x' - x)$ converges to
$x + (1-e)(x'-x) = x' + e(x - x')$, the point
$x_i + (1 - e)(x' - x)$ must be in $U$ for all sufficiently large $i$.
However, none of the points $x_i + (1 - e)(x' - x)$ can belong to $U$
because $f(x_i + (1 - e)(x' - x)) = f(x_i) = y_i$ is not in $f(U)$.
We have a contradiction.
\medskip

(iii) $\Rightarrow$ (ii). Trivial.
\medskip

(ii) $\Rightarrow$ (i). Suppose $Q$ is not LNC and let $x$, $x'$, $x_i$
be a falsifying case. Without loss of generality suppose $x_i + (x' - x)/2$
is never in $Q$ and $x = 0$. Thus $x_i \to 0$ and $x_i + x'/2$ is never in $Q$.
Let $T$ be the projection of $\R^n$ onto $\R^n/[x'] \cong \R^{n-1}$ (the
quotient of $\R^n$ by the one-dimensional subspace spanned by the vector $x'$)
and let $f = T|_Q$.
The next step is to find a relatively open neighborhood $U$ of $x'$ such
that $f(x_i)$ is not eventually in $f(U)$. Then we will have that $f(x_i)$
converges to $f(x') = 0$, but is not eventually in $f(U)$, hence
$f$ is not open.
To fulfill this program let $\phi$ be a linear functional on
$\R^n$ such that $\phi(x') = 1$. Let $U = \{y \in Q : \phi(y) > 3/4\}$.
Then $U$ is a relatively open neighborhood of $x'$ in $Q$. Suppose (by
contradiction) that $f(x_i)$ were eventually in $f(U)$; then there
would exist $z_i \in U$ such that $f(z_i) = f(x_i)$.
By compactness, pass to a subsequence so that $(z_i)$ converges and let
$z = \lim z_i$. Then $\phi(z) =
\lim \phi(z_i) \geq 3/4$. By the definition of $f$ there are real numbers
$t_i$ such that $z_i = x_i + t_i x'$. Then $\phi(z_i) = \phi(x_i) + t_i$. Since
$x_i \to 0$ and $\phi$ is continuous, $\phi(x_i) \to 0$, so eventually
$\phi(x_i) < 1/4$ and thus $t_i > 1/2$. But we have assumed that $x_i + x'/2$
is not in $Q$, so by convexity $x_i + t_i x' = z_i$ cannot be in $Q$ for
$t_i > 1/2$. Since we have chosen $z_i$ to be in $U$, hence in $Q$, we
have a contradiction.\hfill\hal
\bigskip

Next, we consider the simplest situation, where $T$ is the projection of
$\R^n$ onto a quotient by a one-dimensional subspace. Thus,
let $[v]$ be the one-dimensional subspace spanned by a nonzero
vector $v \in \R^n$. Let $T_v: \R^n \to \R^n/[v]$ be the quotient map, let
$Q$ be
a compact, convex subset of $\R^n$, and let $f_v = T_v|_Q: Q \to T_v(Q)$.
In this case there is a simple, explicit
description of a right inverse $g_v: T_v(Q) \to Q$. Namely, given $y = f_v(x)$
with $x \in Q$, let $g_v(y) = x + av$ where
$$a = \inf\{b \in \R: x + bv \in Q\}.$$
Intuitively, if $v$ points ``upward'', then $g_v$ maps $f_v(Q)$ onto the
lower boundary of $Q$. (By replacing $\inf$ with $\sup$, we could similarly
map $f_v(Q)$ onto the upper boundary of $Q$, without essentially affecting
the discussion.)
\medskip

Perhaps surprisingly, the map $g_v$ is not necessarily continuous. Our next
proposition shows that it is continuous for all $v$ precisely if $Q$ is LNC.
\bigskip

\noindent {\bf Proposition 2.} {\it Let $Q$ be a compact, convex subset of
$\R^n$. Then $Q$ is LNC if and only if $g_v: T_v(Q) \to Q$ is continuous
for all nonzero $v \in X$.}
\medskip

\noindent {\it Proof.} Suppose $g_v$ is not continuous for some $v$. Then
there is a sequence $(z_i)$ in $Q$ and an element $z$ in $Q$ such that $y_i =
f_v(z_i)
\to f_v(z) = y$ but $g_v(y_i) \not\to g_v(y)$. Pass to a subsequence so that
$x_i = g_v(y_i)$ converges, say to $x$, and note that $x$ must be of the
form $x = z + av$ since
$$f_v(x) = \lim f_v(x_i) = \lim y_i = f_v(z).$$
Define $x' = g_v(y) = z + a'v$; then $a' < a$ by the minimality of $a'$ in
the definition of $g_v(y)$. Now the LNC property fails for the points $x$
and $x'$ and the sequence $x_i \to x$, because $x_i + (x' - x)/2 =
x_i + (a' - a)v/2$ does not belong to $Q$ for any $i$, by the minimality of
$x_i$. Thus, discontinuity of any $g_v$ implies that $Q$ is not LNC.
\medskip

Conversely, suppose $Q$ is not LNC and find $x, x', x_i \in Q$ such that
$x_i \to x$ but $x_i + (x' - x)/2$ is not in $Q$ for any $i$. Take $v =
x - x'$ and define $y_i = f_v(x_i)$ and $y = f_v(x) = f_v(x')$. Then for
every $i$, $g_v(y_i) = x_i + a_iv$ for some $a_i > -1/2$, since
$x_i + av \not\in Q$ for $a = -1/2$.
However, $g_v(y) = x + av$ for some $a \leq -1$, and so $g_v(y_i) \not\to
g_v(y)$. Thus $g_v$ is not continuous.\hfill\hal
\bigskip

We now give two constructions of continuous sections for linear maps of
LNC sets. The first involves strictly convex norms on $\R^n$.
This means that the closed unit ball is strictly convex. For example,
any euclidean norm on $\R^n$ is strictly convex, as is the $l^p$ norm for
$1 < p < \infty$.
\bigskip

\noindent {\bf Theorem 3.} {\it Let $Q$ be a compact LNC subset
of $\R^n$ and let $T: \R^n \to \R^m$ be a linear map. Let $f = T|_Q: Q \to
T(Q)$.
Fix a strictly convex norm on $\R^n$. Then for each $y \in T(Q)$ there is a
unique element $g(y)$ of $f^{-1}(y) = T^{-1}(y) \cap Q$ with minimal norm,
and the map $g$ is a continuous section of $f$.}
\medskip

\noindent {\it Proof.} There exists an element of $f^{-1}(y)$ with
minimal norm by compactness. For uniqueness, suppose $v$ and $w$ are
distinct and both minimize norm; then $(v + w)/2$ also belongs to the convex
set $f^{-1}(y)$, and it has strictly smaller
norm by strict convexity, a contradiction. So the map $g$ is well-defined.
\medskip

It is also
obviously a section of $f$. To verify continuity, let $(y_i)$ be a sequence
in $T(Q)$ which converges to $y \in T(Q)$ and suppose $(x_i) = (g(y_i))$
fails to converge to $x' = g(y)$. By compactness, we can pass to a subsequence
so that $(x_i)$ converges to some point $x \neq x'$. By LNC, for sufficiently
large $i$ we then have $x_i + (x' - x)/2 \in Q$.
\medskip

However, $f(x) = \lim f(x_i) = y$, and since $x \neq x' = g(y)$ we must have
$\|x\| > \|x'\|$. Therefore
$$\lim \|x_i + (x' - x)/2\| = \|x' + x\|/2 < \|x\| = \lim \|x_i\|,$$
so that $\|x_i + (x' - x)/2\| < \|x_i\|$ for sufficiently large $i$.
Since $f(x' - x) = 0$, this contradicts minimality of the norm of
$x_i = g(y_i)$. This completes the proof.\hfill\hal
\bigskip

We isolate a special case of Theorem 3 in the following corollary; although
simple, we believe this result is new.
\bigskip

\noindent {\bf Corollary 4.} {\it Let $Q$ be a compact, strictly convex subset
of $\R^n$, let $T: \R^n \to \R^m$ be a linear map, and let $f = T|_Q: Q \to
T(Q)$.
For $y \in T(Q)$ define
$g(y)$ to be the unique element of $f^{-1}(y) = T^{-1}(y) \cap Q$ with minimal
euclidean norm. Then $g$ is a continuous section of $f$.}\hfill\hal
\bigskip

If $T$ is a linear map of $\R^n$ into itself (not an essential restriction)
then Theorem 3 and Corollary 4 can be
modified in the following way. Instead of taking $g(y)$ to be the element of
$f^{-1}(y)$ with minimal norm, choose it instead so that $||g(y) - y||$ is
minimized. Trivial modifications of the proofs show that this $g$
is also continuous, if $Q$ is LNC.
\medskip

Our final construction of continuous sections requires some preface. In
previous results we chose a distinguished element of $f^{-1}(y)$ in
simple ways --- in one case by taking one endpoint of a
line segment, and in the other case by minimizing norm. We now introduce
a new, slightly more involved
method of making this choice. This supports our earlier contention
that LNC guarantees that any ``reasonably'' defined section will be
continuous.
\medskip

Let $\{F_r :r = 1, ..., n\}$ be a separating family of linear
functionals on $\R^n$. Thus, for any nonzero $v \in \R^n$ we have $F_r(v)
\neq 0$ for some $r$.
Given a compact set $K$ in $\R^n$, define a nested family of subsets
$K_r$ ($0 \leq r \leq n$) as follows. Let $K_0 = K$. Having defined $K_{r-1}$,
let $K_r$ be the set
$$K_r = K_{r-1} \cap F_r^{-1}(a)$$
where $a = \inf F_r(K_{r-1})$. Geometrically, $K_r$ is
the ``lowest slice'' of $K_{r-1}$, according to $F_r$.
\medskip

Suppose $K_n$ contains two distinct points $v$ and $w$. Then $F_r(v) \neq
F_r(w)$ for some $r$ since the $F_r$'s are separating. So
$v$ and $w$ cannot both belong to $K_r$, a contradiction.
Thus $K_n = \{v\}$ for some $v \in K$. We define $\Gamma(K) = v$; this map
$\Gamma$ performs the desired selection of a distinguished element of $K$.
\bigskip

\noindent {\bf Theorem 5.} {\it Let $Q$ be a compact LNC subset of
$\R^n$, let $T: \R^n \to \R^m$ be a linear map, and let $f = T|_Q: Q \to T(Q)$.
Then $g(y) = \Gamma(f^{-1}(y))$ defines a continuous section of $f$.}
\medskip

\noindent {\it Proof.} It is immediate that $g$ is a section of $f$.
To show that it is continuous, let $y_i \to y$ be a convergent sequence in
$T(Q)$.
Set $x_i = g(y_i)$ and $x' = g(y)$ and suppose $x_i \not\to x'$. Pass to a
subsequence so that $x_i$ converges and let $x$ be its limit. We will show that
$x$, $x'$, and $x_i$ contradict the LNC condition on $Q$.
\medskip

Let $r$ be the smallest index such that $F_r(x) \neq F_r(x')$.
Take $K = f^{-1}(y)$ and observe that both $x$ and $x'$ belong to $K$.
Since $F_s(x) = F_s(x')$ for all
$s < r$ and $K_n = \{x'\}$, it follows that $K_{r-1}$ contains both
$x$ and $x'$. But $K_r$ cannot contain them both, so
$F_r(x) > F_r(x')$.
\medskip

Now for all $i$
set $K^i = f^{-1}(y_i)$, so $K^i_n = \{x_i\}$. Let $v = (x'-x)/2$
and suppose $x_i + v$ is in $Q$. Then $x_i + v$ is in $K^i$. Since
$F_s(v) = 0$ for all $s < r$, it follows that $x_i + v$ is in
$K^i_r$. But $F_r(x_i + v) < F_r(x_i)$, so provided that we
insist on our original assumption that $x \neq x'$ we must then have
$x_i + v \in K^i_r$ and $x_i \not\in K^i_r$,
contradicting the fact that $K^i_n = \{x_i\}$. So we must reject the
assumption that $x_i + v$ is in $Q$.
\medskip

Thus, we have $x_i \to x$ but $x_i + (x' - x)/2 \not\in Q$ for any $i$. This
contradicts the LNC property. We conclude that $x_i$ must have converged to
$x'$, and this shows that $g$ is continuous.\hfill\hal
\bigskip

It is easy to verify that the sections defined in Theorem 5
have the property that $g(y)$ is an extreme point of $f^{-1}(y)$,
for any $y \in T(Q)$. If the range space is $\R^m$ and $Q$ satisfies a
simple geometric condition, this implies that $g(y)$ is actually an extreme
point of $Q$. The condition is that $Q$ should contain no face of dimension
between $1$ and $m$ inclusive; one says that $Q$ has {\it facial dimension}
greater than $m$. (Of course, zero-dimensional faces --- i.e., extreme points
--- cannot be forbidden.) For example, this will be true for all $m < n$ if $Q$
is strictly convex.
\bigskip

\noindent {\bf Theorem 6.} {\it Assume the setup of Theorem 5
and suppose that $Q$ has facial dimension greater
than $m$. Then the continuous section $g$ treated in Theorem 5 takes
values in the extreme points of $Q$. In particular, this will be true
if $Q$ is strictly convex and $m < n$.
\medskip

The same assertion holds for the construction in Theorem 3 and Corollary 4,
provided $0 \not\in Q$.}
\medskip

\noindent {\it Proof.} It is straightforward to verify that $g(y)$ is an
extreme point of $T^{-1}(y) \cap Q$ for any $y \in T(Q)$. Then $g(y)$ is
an extreme point of $Q$ by [2, Theorem 1.6].\hfill\hal
\bigskip
\bigskip

\noindent {\bf 3. Finite-dimensional examples}
\bigskip

We now list several elementary finite-dimensional examples and counterexamples.
\bigskip

\noindent {\bf A. Strictly convex sets.} It is immediate from the definitions
that every strictly convex set is LNC. (This is true in infinite dimensions
as well.)
\bigskip

\noindent {\bf B. Sets in $\R^1$ and $\R^2$.} Every convex set in
$\R^1$ or $\R^2$ is LNC.
\bigskip

\noindent {\bf C. Polytopes.} The following proposition is trivial.
\bigskip

\noindent {\bf Proposition 7.} {\it Let $Q$ and $Q'$ be LNC sets.
Then $Q \cap Q'$ is also LNC.}\hfill\hal
\bigskip

\noindent As a consequence of this result and the easy observation that
any half-space is LNC, it follows that convex polytopes in $\R^n$ are LNC.
This provides us with a large class of LNC sets which are not strictly
convex. Intersecting polytopes with strictly convex sets provides even more
examples.
\medskip

Note that the intersection of infinitely many LNC sets need not be LNC.
Indeed, any closed convex set is an intersection of half-spaces, so this
follows just from the fact that there exist closed convex sets which are
not LNC. We give examples of such sets in Sections E and F below.
\bigskip

\noindent {\bf D. Images of compact LNC sets.} 
Any linear image of a compact LNC set in $\R^n$ is also LNC.
\bigskip

\noindent {\bf Proposition 8.} {\it Let $Q$ be a compact LNC set in $\R^n$
and let $T: \R^n \to \R^m$ be linear. Then $T(Q)$ is also a compact LNC set.}
\medskip

\noindent {\it Proof.} $T(Q)$ is compact because $T$ is continuous. To verify
that $T(Q)$ is LNC, let $y, y', y_n \in T(Q)$ and suppose $y_n \to y$. Find
$x', x_n \in Q$ such that $T(x') = y'$ and $T(x_n) = y_n$. By passing to a
subsequence we may assume that $(x_n)$ converges; letting $x = \lim x_n$,
we have $T(x) = \lim T(x_n) = y$. Now $x_n + (x' - x)/2$ is eventually in $Q$
because $Q$ is LNC, and applying $T$ shows that $y_n + (y' - y)/2$ is
eventually in $T(Q)$.\hfill\hal
\bigskip

\noindent {\bf E. Cones.} 
\bigskip

\noindent {\bf Example 9.}
The simplest example of a compact, convex, non-LNC set in $\R^n$ is a right
circular cone, explicitly given (for example) as the convex hull of the set
$\{(1 - \cos t, \sin t, 1): 0 \leq t \leq 2\pi\}$ together with the origin
in $\R^3$. The LNC condition is falsified by the points $z = (0, 0, 1)$,
$z' = (0, 0, 0)$, and $z_n = (1 - \cos(1/n), \sin(1/n), 1)$. Taking
$T: \R^3 \to \R^2$ to be the projection $T(x,y,z) = (x,y)$, the corresponding
section $g$ defined by any of the constructions in Section 2 is discontinuous.
In every case, $g(0, 0) = (0, 0, 0)$ while
$g(1-\cos t, \sin t) = (1 - \cos t, \sin t, 1)$ for $0 < t < 2\pi$.
\bigskip

It is worth noting that the restriction of $T$ to the cone does have
continuous right inverses, however; the simplest is the map $(x,y) \mapsto
(x,y,1)$. A closer analysis shows that in fact the restriction to the cone
of any linear map $T$ from $\R^3$ to $\R^2$ (indeed, to any $\R^m$) has a
continuous section. The following is an example where no continuous sections exist.
\bigskip

\noindent {\bf Example 10.} Let
$Q$ be the convex hull of the helix $\{(\cos t, \sin t, t): 0 \leq t \leq 2\pi\}$
in $\R^3$ and consider the orthogonal projection $T: \R^3 \to \R^2$ onto
the $xy$-plane. Then any right inverse $g$ of $f= T|_Q$ must satisfy
$g(\cos t, \sin t) = (\cos t, \sin t, t)$ for $0 < t < 2\pi$, and hence
must have a discontinuity at $(1,0)$.
\bigskip

Example 9 has the following generalization.
\bigskip

\noindent {\bf Proposition 11.} {\it Let $Q_0$ be a compact, convex subset
of $\R^n$ and let
$$Q = \{(tx, t): x \in Q_0, 0 \leq t \leq 1\}$$
be the suspension of $Q_0$ in $\R^{n+1} = \R^n \times \R$. Then $Q$ is LNC
if and only if $Q_0$ is a polytope.}
\medskip

\noindent {\it Proof.} If $Q_0$ is a polytope then so is $Q$, hence $Q$ is
LNC. Conversely, if $Q_0$ is not a polytope then it has infinitely many
extreme points. (Any compact, convex set is the closed hull of its extreme points,
so if it had only finitely many extreme points it would be a polytope.) Let
$(x_i)$ be a sequence of distinct
extreme points of $Q_0$ which converges to a cluster
point $x$. Then in $Q$ we have $(x_i, 1) \to (x, 1)$; also $(0, 0)$ is in $Q$,
so the LNC condition would require that $(x_i - x/2, 1/2)$ is eventually in
$Q$. But this cannot be, because
if $(y/2, 1/2) = (x_i - x/2, 1/2)$ is in $Q$ then $y = 2x_i - x$ is in $Q_0$.
Then $x_i = (y + x)/2$, contradicting the fact that $x_i$ is extreme.\hfill\hal
\bigskip

\noindent {\bf F. Matrix algebras.} The unit balls of the most common
finite-dimensional Banach spaces are either strictly convex or polytopes,
and hence are LNC by reasons given above. However, finite-dimensional unit
balls need not be LNC.
\bigskip

\noindent {\bf Proposition 12.} {\it Neither the unit ball nor the
positive part of the unit ball of the $n\times n$ matrix algebra (with
operator norm or trace norm) is LNC for $n \geq 2$.}
\medskip

\noindent {\it Proof.} Take $n = 2$. To show that the positive part of the
unit ball is not LNC, define $x'$ to be $0$, $x = \pmatrix{1&0\cr 0&0\cr}$,
and, for $t \in (0,1)$, define
$$x_t = \pmatrix{t& \sqrt{t - t^2}\cr \sqrt{t - t^2}& 1-t\cr}.$$
Then $x_t \to x$ as $t \to 1$, but for any $t \in (0,1)$ the matrix
$x_t + (x' - x)/2$ has negative determinant, and hence does not belong to the
positive part of the unit ball.
\medskip

Replacing each matrix $A$ with $I - A$ in the above argument shows that
the unit ball of the $2\times 2$ matrix algebra is not LNC, and the same
construction can be carried out in the upper left $2\times 2$ block of any
$n \times n$ matrix algebra.\hfill\hal
\bigskip

\noindent {\bf G. Sets without one-dimensional faces.} The failures of the
LNC condition in the cone and matrix examples considered above all happen on
one-dimensional faces of $Q$. However, this is not essential. We now
describe a set in $\R^4$ which has no one-dimensional faces but still fails to
be LNC.
\medskip

\noindent {\bf Example 13.}
Let $Q_0$ be a cone in $\R^3$ and let $I$ be a compact interval. Then
$Q_0 \times I$, the Cartesian product which is contained in $\R^4$, is not
LNC. In fact, for any $t \in I$ the intersection of $Q_0 \times I$ with
$\R^3 \times \{t\}$ is the ``slice'' $Q_0 \times \{t\}$, which is
isometric to $Q_0$, and hence is not LNC. By Proposition 7 we
conclude that $Q_0 \times I$ cannot be LNC.
\medskip

Now $Q_0 \times I$ has some one-dimensional faces: for any one-dimensional
face $F$ of $Q_0$, the sets $F \times \{a\}$ and $F \times \{b\}$ are
one-dimensional faces of $Q_0 \times I$, where $I = [a, b]$. Also, for any
extreme point $x$ of $Q_0$, the set $\{x\} \times I$ is a one-dimensional
face of $Q_0 \times I$. However, all of these faces can be ``sliced off'' to
get a set with the desired properties.
\medskip

For concreteness, let $Q_0$ be the convex hull of the point $(0, 0, 10)$
and the circle $(\cos t, \sin t, -10)$ and let $I$ be the interval
$[-10, 10]$. Let $Q$ be the intersection of $Q_0 \times I$ with the 4-ball
of radius 2 about the origin. This removes all of the one-dimensional faces
of $Q_0 \times I$ --- the sphere of radius 2 does not contain any lines, so
there are no one-dimensional faces on the boundary, and inside the sphere
$Q$ is a part of $Q_0 \times I$ which has only two-dimensional faces.
Furthermore,
the intersection of $Q$ with $\R^3 \times \{0\}$ is a truncation of $Q_0$
which is not LNC. So as before, it follows that $Q$ is not LNC.
\bigskip
\bigskip

\noindent {\bf 4. Closed LNC sets}
\bigskip

The results in Section 2 all trivially have infinite-dimensional analogs
in any TVS
(just replace sequences by nets). However, this is not very interesting
because infinite-dimensional compact LNC sets do not exist (Theorem
23). Thus, the consequences for compact sets in infinite dimensions are
negative: there always exist linear maps whose restrictions are not open,
and sections defined in various ways are not continuous in general.
\medskip

But LNC sets which are merely closed are easy to construct in infinite
dimensions (see the first paragraph of Section 5). And while the proofs
are more difficult, we do have analogs of the results in Section 2, provided
that $T$ is left Fredholm, i.e.\ $T$ is a bounded operator with
closed range and finite-dimensional kernel.
Of course, if the TVS $X$ is already finite-dimensional then this is no restriction
at all. Thus one could say that in finite dimensions closed LNC sets are
practically as well-behaved as compact LNC sets.
\medskip

Our fundamental tool is Lemma 15. In its proof we need to use Hilbert
space techniques. This is irrelevant to the finite-dimensional
case because without loss of generality we can always equip $\R^n$
with a Euclidean norm. We believe that in infinite dimensions, the
Hilbert space condition on $X$ can be weakened to $X$ being a uniformly
convex Banach space. This issue will be addressed in a future paper.
\medskip

It is worth noting that if the closed LNC set in question is bounded, then
the results in this section are easy variations on the results in Section 2.
The idea is that if $(x_n)$ is a bounded sequence in a Hilbert space and
$(Px_n)$ converges, where $P$ is an orthogonal projection with finite-dimensional
kernel, then we can pass to a subsequence so that $(P^\perp x_n)$ also converges,
and then $(x_n)$ converges because $x_n = Px_n + P^\perp x_n$. This technique is
used in the proof of Theorem 20.
\medskip

Our proof of Lemma 15 requires the following fundamental result from
[8] and [9]. For the reader's convenience we include an easy proof here.
\bigskip

\noindent {\bf Theorem 14.} {\it Let $Q$ be an LNC set in a TVS and let
$x,y \in Q$ be distinct. Then $p = {1\over 2}(x + y)$ has a relative
neighborhood $U$ in $Q$ such that any $q \in U$ lies in the interior of
a line segment which is contained in $Q$ and parallel to $[x,y]$.}
\medskip

\noindent {\it Proof.} Suppose the conclusion fails. Then there is a net
$(q_t) \subset Q$ which converges to $p$ such that no $q_t$ is
in the interior of a line segment which is contained in $Q$ and parallel to
$[x,y]$. Thus, setting $v = {1\over 4}(y - x)$, for each $t$ the points
$q_t - v$ and $q_t + v$ cannot both belong to $Q$. By passing to
a subnet, without loss of generality we may suppose that $q_t + v
\not\in Q$ for all $t$. Then the substitutions $x = p$, $x_t =
q_t$, and $y = y$ contradict the fact that $Q$ is LNC.\hfill\hal
\bigskip

In any Hilbert space (and hence in any finite-dimensional space) the
assertion of Theorem 14 is actually equivalent to the LNC property
[9, Theorem 2.2].
\medskip

We use the notation $[X]_\epsilon$ for the closed $\epsilon$-ball
about the origin in $X$. We use the term
``relative interior'' to mean the interior of a convex
set relative to the affine subspace of $X$ that it spans.
\bigskip

\noindent {\bf Lemma 15.} {\it Let $X$ be a Hilbert space, let $Q \subset X$
be a closed LNC set, and let $P: X \to X$ be an orthogonal projection with
finite-dimensional kernel. Suppose $0 \in Q$ and $0$ is in the relative
interior of $Q \cap {\rm ker}\, P$. Let $(x_i) \subset Q$ and suppose
$P(x_i) \to 0$. Then $P'(x_i) \to 0$, $(P'(x_i))$ is eventually in $Q$, and
$P' \geq P$, where $P'$ is the orthogonal projection onto
$({\rm span}\, (Q\cap {\rm ker}\, P))^\perp$.}
\medskip

\noindent {\it Proof.}  Let $W = {\rm span}\, (Q\cap {\rm ker}\, P)$. By
hypothesis, for every $v \in W$ there is a line segment in the direction
of $v$ which is contained in $W$ and
contains $0$ in its interior. Fixing an orthonormal basis
$\{v_1, \ldots, v_n\}$ of $W$, it follows from $n$ applications of Theorem 14
that we can choose $\epsilon' > 0$ and $\delta' > 0$ such that for any
$q \in Q \cap [X]_{\epsilon'}$
we have $[q + \delta' v_i, q - \delta' v_i] \subset Q$ for $1 \leq i \leq n$.
Letting $\delta = \delta'/({\rm dim}\, W)^{1/2}$, the ball of
radius $\delta$ about $q$ is contained in the convex hull of the line segments
$[q + \delta' v_i, q - \delta' v_i]$, so
this implies that for every $q \in Q \cap [X]_{\epsilon'}$
we have $q + [W]_\delta \subset Q$. Let $\epsilon = {\rm min}(\epsilon', 1)$.
\medskip

Let $Q' = Q \cap [X]_\epsilon$ and define $H$ to be the intersection of all
closed half-spaces $H_\beta$ in $X$ such that $Q \subset H_\beta$ and
$Q' \cap \partial H_\beta \neq \emptyset$. If $p \not\in Q$ and $\|p\| \leq
\epsilon$, we claim that there exists an $H_\beta$ which excludes $p$.
To see this, let $q$ be the unique element of $Q$ such that $\|q - p\|$ is
minimized [7, Theorem 12.3]. Notice that if $\|q\| > \|p\|$ then
$q' = (\langle p,q\rangle/\|q\|^2)q$ is in $Q'$ (since
$\|q'\| \leq \|p\| \leq \epsilon$, hence $\langle p,q\rangle/\|q\|^2 \leq
\epsilon$, and $0,q \in Q$). Also $\|q' - p\| < \|q - p\|$ since
$q'$ is the projection of $p$ onto the one-dimensional subspace $[q]$;
this contradicts the minimality of $\|q - p\|$, and so we must have
$\|q\| \leq \|p\|$. Thus $q \in Q'$. Now define
$$\eqalign{H_\beta
&= \{x \in X: \langle x-p, q-p\rangle \geq \langle q-p, q-p\rangle\}\cr
&= \{x \in X: \langle x, q-p\rangle \geq \langle q, q-p\rangle\}.\cr}$$
We have $Q \subset H_\beta$. To see this, let $x \in Q$. Since $q$ is
also in $Q$, so is $tx + (1-t)q$ for any t in [0,1]. Since $q$ minimizes
distance to $p$, we therefore have
$$\|t(x-q)  +  q  -  p\| = \|tx + (1-t)q  -  p\|   \geq   \|q - p\|$$
for all t in [0,1].  Define $f(t) =  \|t(x-q) + q - p\|^2$. By the above
we must have $f'(0) \geq 0$.  But $f'(t) = 2t \|x-q\|^2 +
2 \langle x - q, q - p\rangle$. So $f'(0) = 2 \langle x - q, q - p\rangle
\geq 0$.  That is, $\langle x, q - p\rangle \geq \langle q, q - p\rangle$.
Thus $x$ is in $H_\beta$, and we have shown $Q \subset H_\beta$.
It is straightforward to verify that $q \in
Q' \cap \partial H_\beta$ and $p \not\in H_\beta$. This proves the claim.
\medskip

Thus $Q' = H\cap [X]_\epsilon$. Now each $\partial H_\beta$ contains an
element $q$ of $Q'$, and such an element is in the relative interior of
$Q \cap (q + W)$. Therefore each $\partial H_\beta$ contains a translate
of $W$ and from this it follows that $H = H + W$.
Also, $H \cap {\rm ker}\, P = W$: for any $p \in {\rm ker}\, P$, $p \not\in W$,
$\|p\| \leq {\rm min}(\epsilon, \delta)$, we must have $p \not\in Q$ (since $p \not\in W$),
therefore $p \not\in Q'$, and this implies that $p \not\in H$ because
$Q' = H\cap [X]_\epsilon$. This
proves that $H \cap {\rm ker}\, P \subset W$; the reverse inclusion is trivial.
\medskip

Now let $P'$ be the orthogonal projection of $X$ onto $W^\perp$. Note immediately
that $P' \geq P$, that is, $P\circ P' = P' \circ P = P$ because $W \subset
{\rm ker}\, P$. If $x_i \to 0$ we are done, because then $\|P'(x_i)\| \leq
\|x_i\| \to 0$ and $\|x_i\| \leq {\rm min}(\epsilon, \delta)$ implies that 
$x_i \in Q'$ and $\|x_i - P'(x_i)\| \leq \delta$, hence $P'(x_i) \in Q$.
Otherwise, pass to a subsequence so that $\|x_i\| \geq \gamma > 0$ for all $i$.
Let $y_i = (P' - P)(\gamma_i x_i)$ where $\gamma_i = {\rm min}(\epsilon/\|x_i\|, 1)$
(and thus $\gamma_i \leq \epsilon/\gamma$). Observe that $P(y_i) = 0$ and $P'(y_i) = y_i$.
Pass to another subsequence to ensure that
$(y_i)$ converges to some point $y$ (this can be done because $\|y_i\| \leq \epsilon$
for all $i$ and ${\rm ker}\, P$ is finite-dimensional) and note
that $P'(\gamma_i x_i) \to y$ since $P(\gamma_i x_i) \to 0$. Since
$\gamma_i x_i \in Q' \subset H$
for all $i$, we also have $P'(\gamma_i x_i) \in H + W = H$.
It follows that $y \in H$ as well. Also, $P(y_i) = 0$ implies that $P(y) = 0$,
so by the result of the last paragraph we have
$y \in H \cap {\rm ker}\, P = W$. Yet $P'(y_i) = y_i$ implies $P'(y) = y$, so
$y \in W^\perp$, and therefore $y = 0$. Therefore $P'(\gamma_ix_i) \to 0$,
and this implies that $P'(x_i) \to 0$.
\medskip

Finally, since $P'(x_i) \to 0$ this sequence eventually lies in $[X]_\epsilon$.
It also belongs to $H$ because $x_i \in Q \subset H$ and $H + W = H$. Since
$Q' = H \cap [X]_\epsilon$, we conclude that $P'(x_i) \in Q'$ eventually.
In particular, $P'(x_i)$ is eventually in $Q$.\hfill\hal
\bigskip

Now we proceed to the promised analogs of the results in Section 2.
\bigskip

\noindent {\bf Theorem 16.} {\it Let $Q$ be a closed, convex subset of
a Hilbert space $X$. The following conditions are equivalent:
\medskip

{\narrower{
\noindent (i) $Q$ is LNC;
\medskip

\noindent (ii) for every Hilbert space $Y$ and every left Fredholm map
$T: X \to Y$, the restriction map $f = T|_Q: Q \to T(Q)$ is open.
\medskip}}}

\noindent {\it Proof.} (i) $\Rightarrow$ (ii). Suppose $Q$ is LNC and let
$T: X \to Y$ be left Fredholm.
Let $P$ be the orthogonal projection with ${\rm ker}\, P = {\rm ker}\, T$.
Then $T = T'\circ P$ where $T'$ is a linear homeomorphism from
${\rm ran}\, P$ onto ${\rm ran}\, T$. So we only need to show that
$P|_Q: Q \to P(Q)$ is open.
\medskip

Suppose $P|_Q$ is not open; then there exist $x \in Q$, $(x_i) \subset Q$,
and a relatively open set $U \subset Q$ containing $x$ such that $P(x_i) \to
P(x)$ but $P(x_i) \not\in P(U)$ for all $i$. Let $V = {\rm ker}\, P$; then
$(x + V) \cap Q$ is a closed convex set, and there exist points in its
relative interior arbitrarily close to $x$. Thus we can find a point $x'$
belonging to the intersection of $U$ with the relative interior of
$(x + V) \cap Q$. Note that $x' - x \in V$, so $P(x') = P(x)$. Thus,
replacing $x$ with $x'$, we may assume that $x$ belongs to the relative
interior of $(x + V) \cap Q$. Translating $Q$ by $x$, we can further assume
that $x = 0$. The hypotheses of Lemma 15 are now satisfied, so we have
$P'(x_i) \to 0$ and $P'(x_i) \in Q$ eventually. Therefore $P'(x_i) \in U$
eventually, so that $P(x_i) = P(P'(x_i)) \in P(U)$ eventually, contradicting
the choice of $U$. So $P|_Q$ must be open.
\medskip

\noindent (ii) $\Rightarrow$ (i). The proof of Theorem 1 (ii) $\Rightarrow$ (i)
works for any Hilbert space $X$ in place of $\R^n$, modulo one minor modification.
Instead of setting $U = \{y \in Q: \phi(y) > 3/4\}$, define
$U = \{y \in Q: 5/4 > \phi(y) > 3/4\}$. Then regardless of whether $Q$ is
compact, we can pass to a subsequence of $(z_i)$ so that $\phi(z_i)$ converges,
and this is sufficient to complete the proof.\hfill\hal
\bigskip

For any nonzero $v \in X$, let $T_v$ be the orthogonal projection onto the
orthocomplement of $v$. Let $Q$ be a closed, convex subset of $X$ and
let $f_v = T_v|_Q: Q \to T_v(Q)$. For any $x \in T_v(Q)$, the set
$\{b \in \R: x + bv \in Q\}$ is a closed interval in $\R$ and so, although
it may not have a smallest element, it does contain
a unique element of minimal absolute value. Let $a$ be this number, and
define $g_v(x) = x + av$. This is a slightly modified version of the map
$g_v$ defined in Section 2 which is necessary if $Q$ is closed and unbounded.
The original definition would work if $Q$ were assumed to be bounded, but
it is ill-defined in general
because $\inf\{b \in \R: x + bv \in Q\}$ may not exist. The disadvantage
of our new definition of $g_v$ is that if $Q$ is translated in $X$ the
corresponding $g_v$ may not be a translate of the original $g_v$.
\bigskip

\noindent {\bf Proposition 17.} {\it Let $Q$ be a closed, convex subset of
a Hilbert space $X$. Then $Q$ is LNC if and only if $g_v: T_v(Q') \to Q'$ is
continuous for all nonzero $v \in X$ and all translates $Q'$ of $Q$, where
$g_v$ and $T_v$ are as defined above.}
\medskip

\noindent {\it Proof.} We can assume $Q' = Q$ in the forward direction of
the proof. Thus, suppose $Q$ is LNC and $g_v$ is not continuous for some
$v \in X$. Then there is a sequence $(y_i)$ in $T_v(Q)$ and an element $z$ in
$Q$ such that $y_i \to f_v(z) = y$ but $g_v(y_i) \not\to g_v(y)$.
We may assume that $z$ is in the relative interior of $(z + {\rm ker}\, T_v) \cap Q$.
\medskip

Let $x_i = g_v(y_i)$. Applying Lemma 15 to the sequence $(x_i - z)$, the LNC
set $Q - z$, and the projection $P = T_v$,
we conclude that $P'(x_i) \to P'(z)$ and $P'(x_i) \in Q$ eventually. But since
${\rm ker}\, T_v$ is one-dimensional and $P' \geq P$, either $P'$ is the identity operator
or $P' = P$. In the latter case, $y_i = P'(x_i) \in Q$ eventually, so that
$g_v(y_i) = y_i$ eventually and $y \in Q$, hence $g_v(y) = y$, and therefore
$g_v(y_i) \to g_v(y)$, contradicting the choice of $(y_i)$ and $y$.
So we can assume that $P'$ is the identity operator and so $x_i \to z$.
Since $P'$ is the identity, it follows that $(z + {\rm ker}\, T_v) \cap Q =
\{z\}$, so we must have $g_v(y) = z$. Therefore $g_v(y_i) = x_i \to z =
g_v(y)$, again contradicting the choice of $(y_i)$ and $y$. We conclude that
if $Q$ is LNC then $g_v$ is continuous.
\medskip

Conversely, suppose $Q$ is not LNC and find $x, x', x_i \in Q$ such that
$x_i \to x$ but $x_i + (x' - x)/2$ is not in $Q$ for any $i$. Let $Q' = Q - x'$
and define $z = v = x - x'$, $z' = 0$, and $z_i = x_i - x'$.
Then define $y_i = f_v(z_i)$ and $y = 0$. Since $z_i \to z = v$, it follows that
$y_i = T_v(z_i) \to 0 = y$. Also $g_v(y) = 0$, whereas $g_v(y_i) = z_i + a_iv$
where $a_i > -1/2$ since $z_i \in Q'$ but $z_i -v/2 = x_i + (x' - x)/2 - x'
\not\in Q'$. But convergence of $g_v(y_i)$ to $0 = z - v$, together with the
fact that $z_i \to z$, would imply that $a_i \to -1$, a contradiction. Thus
$g_v(y_i) \not\to g_v(y)$, and so $g_v$ is not continuous.\hfill\hal
\bigskip

\noindent {\bf Theorem 18.} {\it Let $Q$ be a closed LNC subset of a Hilbert
space $X$ and let $T$ be a left Fredholm map from $X$ into another Hilbert
space $Y$.
Let $f = T|_Q: Q \to T(Q)$. Then for each $y \in T(Q)$ there is a
unique element $g(y)$ of $f^{-1}(y) = T^{-1}(y) \cap Q$ with minimal norm,
and the map $g$ defined in this way is a continuous section of $f$.}
\medskip

\noindent {\it Proof.} As in the proof of Theorem 16 (i) $\Rightarrow$ (ii),
we may assume that $T$ is the orthogonal projection $P$ of $X$ onto a subspace.
\medskip

Existence and uniqueness of $g(y)$ is Theorem 12.3 of [7].
To verify continuity, let $(y_i)$ be a sequence in $T(Q)$ which
converges to $y \in T(Q)$ and suppose $(x_i) = (g(y_i))$ fails to converge
to $x' = g(y)$. Fix $\epsilon > 0$ and find $x$ in the relative
interior of $(x' + {\rm ker}\, T) \cap Q$ such that $\|x' - x\| \leq \epsilon$.
Then apply Lemma 15 to $Q' = Q - x$ and the sequence $(x_i') = (x_i - x)$; we get
$P'(x_i) \to P'(x)$ and $P'(x_i - x) \in Q'$ eventually, hence
$w_i = P'(x_i) + x - P'(x) \in Q$ eventually. Note that
$w_i \to x$ and $T(w_i) = y_i$, so $\|w_i\| \geq \|x_i\|$ eventually. Since
$(w_i)$ is convergent and hence bounded, we can pass to a subsequence to ensure
that the norms $\|x_i\|$ converge. Then
$$\|x'\| + \epsilon \geq \|x\| = \lim \|w_i\| \geq \lim \|x_i\|.$$
But $Q$ is closed and convex, hence it is weakly closed, so there is a weak
cluster point $x'' \in Q$ of $(x_i)$; and
$$T(x'') = \lim T(x_i) = \lim y_i = y.$$
So we also have
$$\|x'\| \leq \|x''\| \leq \lim \|x_i\|,$$
which together with the above (and the fact that $\epsilon$ is arbitrary)
establishes that $\|x'\| = \lim \|x_i\|.$
\medskip

Again fix $\epsilon > 0$ and, as in the last paragraph, find a sequence $(w_i)$
which is eventually in $Q$, satisfies $T(w_i) = y_i$, and converges to $x \in Q$ 
where $T(x) = y$ and $\|x' - x\| \leq \epsilon$. Pass to a subsequence so that
$w_i \in Q$, $\|w_i - x\| \leq \epsilon$, and $|\, \|x_i\| - \|x'\|\, | \leq
\epsilon$ for all $i$. Then $\|w_i\| \geq \|x_i\|$ and
$$\|w_i\| \leq \|x\| + \epsilon \leq \|x'\| + 2\epsilon
\leq \|x_i\| + 3\epsilon$$
for all $i$. Also, since both $w_i$ and $x_i$ are in the convex set
$Q \cap T^{-1}(y_i)$, and $x_i$ is the unique element of this set with minimal
norm, it follows that $\langle w_i, x_i\rangle \geq \|x_i\|^2$. Therefore
$$\eqalign{\|w_i - x_i\|^2 & = \|w_i\|^2 -2\langle w_i, x_i\rangle + \|x_i\|^2\cr
& \leq (\|x_i\| + 3\epsilon)^2 - 2\|x_i\|^2 + \|x_i\|^2\cr
& = 3\epsilon(2\|x_i\| + 3\epsilon)\cr
& \leq 3\epsilon(2\|x'\| + 5\epsilon)\cr}$$
eventually.
Thus, by choosing $\epsilon$ sufficiently small we can ensure that $x_i$ and
$w_i$ are eventually arbitrarily close and simultaneously
that $x = \lim w_i$ is arbitrarily
close to $x'$. This implies that $(x_i)$ converges to $x'$. Since we passed
to a subsequence, we have really shown that every subsequence of
the original sequence $(x_i)$ has a subsequence which converges to $x'$. But this
implies that the original sequence $(x_i)$ converges to $x'$, contradicting our
assumption that $g$ is not continuous. Therefore $g$ is continuous.\hfill\hal
\bigskip

Let $T: X \to Y$ be a left Fredholm map.
Let $V = {\rm ker}\, T$ and choose a separating family of linear functionals
$\{F_r: r = 1, \ldots, n\}$ on $V$. Given a compact set $K$ in $V$, define
the sets $K_r$ ($0 \leq r \leq n$) and the element $\Gamma(K)$ as in Section 2,
just before Theorem 5. For any parallel affine subspace $V'$, write $V' = V + v$ with
$v \perp V$, and for $K \subset V'$ compact define $K_r = (K - v)_r + v$ and
$\Gamma(K) = \Gamma(K - v) + v$.
\medskip

(Actually, the functionals $F_r$ can be chosen independently of $T$; let
$\{F_r\}$ be any well-ordered separating set of linear functionals on $X$ and
use them to construct $\Gamma(K)$. In effect we are using the finite set
$\{F_{k_i}\}$, $1 \leq i \leq {\rm dim}\, V$, where $k_1 = 1$ and $k_{i+1}$
is the first index such that $\{F_{k_1}|_V, \ldots, F_{k_{i+1}}|_V\}$ is
linearly independent.)
\medskip

Before proving the next result, we give a counterexample which shows that
in contrast to previous results in this section, using
$\Gamma$ to define a continuous section of $T|_Q$ requires that $Q$ be bounded.
This is true even in finite dimensions. Of course, there is an immediate problem in the
construction of $K_r$ if $K$ is unbounded because then there may be no ``lowest
slice''; even if $a = \inf F_r(K_{r-1})$ is finite, $K_{r-1} \cap F_r^{-1}(a)$
may be empty. Moreover, even if we can define a section using this procedure it
need not be continuous.
\bigskip

\noindent {\bf Example 19.} Define
$$Q = \{(x,y,z) \in \R^3: x \geq 0, y \geq 0, x + y \leq 1,{\hbox{ and }}
z \geq {{(1 - y)^3}\over x}\},$$
where we set $(1 - y)^3/x = 0$ when $x = 0$ and $y = 1$.
Define $T: \R^3 \to \R$ by $T(x,y,z) = x$. One can verify that the lower
boundary of $Q$ is strictly convex, by checking that the Hessian of the function
$f(x,y) = (1-y)^3/x$ is strictly positive-definite. It follows from this that $Q$
is LNC. But $T(Q) = [0,1]$, and if we take $F_1(y,z) = y$ and $F_2(y,z) = z$ then
we have $\Gamma(f^{-1}(x)) = (x,0, 1/x)$
for $x \in (0,1]$ and $\Gamma(f^{-1}(0)) = (0,1,0)$. Thus $\Gamma$ is not continuous.
\bigskip

\noindent {\bf Theorem 20.} {\it Let $Q$ be a closed, bounded LNC subset of
a Hilbert space $X$, let $T: X \to Y$ be left Fredholm, and let $f = T|_Q: Q \to T(Q)$.
Then $g(y) = \Gamma(f^{-1}(y))$ defines a continuous section of $f$.}
\medskip

\noindent {\it Proof.} Observe that $f^{-1}(y) = T^{-1}(y) \cap Q$ is the
intersection of a finite-dimensional affine subspace of $X$ with a closed, bounded set,
so it is compact, and therefore $\Gamma(f^{-1}(y))$ is defined for all $y \in T(Q)$.
As usual, we may assume that $T = P$ is an orthogonal projection with finite-dimensional
kernel. Now let $(y_i) \to y$ be a convergent sequence in $T(Q)$.
Set $x_i = g(y_i)$ and $x' = g(y)$. Since
$x_i = P(x_i) + P^\perp(x_i)$ and ${\rm ker}\, P = {\rm ran}\, P^\perp$ is
finite-dimensional, we can pass to a subsequence so that $P^\perp(x_i)$
converges; then $x_i$ also converges, say to $x$. Note that $T(x) =
\lim T(x_i) = y = T(x')$. We will show that
$x$, $x'$, and $x_i$ contradict the LNC condition on $Q$ unless $x = x'$,
and therefore $g$ must be continuous.
\medskip

Thus suppose $x \neq x'$, and
let $r$ be the smallest index such that $F_r(x) \neq F_r(x')$.
Take $K = f^{-1}(y)$. Since $F_s(x) = F_s(x')$ for all
$s < r$ and $K_n = \{x'\}$, it follows that $K_{r-1}$ contains both
$x$ and $x'$. But $K_r$ cannot contain them both, so
$F_r(x) > F_r(x')$.
\medskip

Now for all $i$
set $K^i = f^{-1}(y_i)$, so $K^i_n = \{x_i\}$. Let $v = (x'-x)/2$
and suppose $x_i + v$ is in $Q$. Then $x_i + v$ is in $K^i$. Since
$F_s(v) = 0$ for all $s < r$, it follows that $x_i + v$ is in
$K^i_r$. But $F_r(x_i + v) < F_r(x_i)$, so provided that we
insist on our original assumption that $x \neq x'$ we must then have
$x_i + v \in K^i_r$ and $x_i \not\in K^i_r$,
contradicting the fact that $K^i_n = \{x_i\}$. So we must reject the
assumption that $x_i + v$ is in $Q$.
\medskip

Thus, we have $x_i \to x$ but $x_i + (x' - x)/2 \not\in Q$ for any $i$. This
contradicts the LNC property. We conclude that $x_i$ must have converged to
$x'$, and this shows that $g$ is continuous.\hfill\hal
\bigskip

Throughout this section we have insisted that $T$ have closed range and
finite-dimensional kernel. Are these assumptions necessary? The first
certainly is: let $Q$ be the unit ball of $l^2 = l^2(\N)$ and define
$T: l^2 \to l^2$ by $T((a_n)) = (a_n/(n+1))$. Then $T$ has null kernel but
does not have closed range. There is only one right inverse $g: T(Q) \to Q$
of $f = T|_Q$, and it is not continuous.
\medskip

We believe that finite-dimensionality of ${\rm ker}\, T$ is also necessary
in general,
but we have not found a counterexample. The most extreme counterexample would
involve a bounded, strictly convex set $Q$ and a map $T$ with cofinite-dimensional
kernel, but even this case remains open. We pose it as a problem.
\bigskip

\noindent {\bf Open problem.} Let $Q$ be a closed, bounded, strictly convex
subset of a Hilbert space $X$ and let $T: X \to \R^n$ be a bounded linear map.
Is $T|_Q: Q \to T(Q)$ necessarily open?
\bigskip
\bigskip

\noindent {\bf 5. Infinite-dimensional counterexamples.}
\bigskip

Infinite-dimensional LNC sets are abundant.
As in the finite-dimensional case, strictly convex sets and half-spaces are
always LNC, and so are finite intersections and Cartesian products of these
sets. However, it is difficult to find other examples than these (the
unit ball of $c_0$ is an example [9, Theorem 1.1]).
\bigskip

\noindent {\bf A. Compact sets in infinite dimensions.}
\bigskip

We now prove that no infinite-dimensional compact convex set is LNC.  We wish
to thank Stephen Simons for essential contributions to our original proof
of this result.
\bigskip

\noindent {\bf Lemma 21.} {\it Let $Q$ be the unit ball of any                            
infinite-dimensional dual Banach space $X^*$. With respect to the weak*
topology, $Q$ is not LNC.}
\medskip
 
\noindent {\it Proof.} First, find a sequence of elements $x_n \in X$ and
$f_n \in X^*$ with the properties $\|x_n\| = \|f_n\| = 1$, $f_n(x_n) = 1$,
and $f_n(x_k) = 0$ if $k < n$. This can be done inductively. The base step
$n = 1$ is trivial. For $n = k + 1$ let $Y$ be the span of $x_1, \ldots, x_k$;
find $x_{k+1} \in X$ such that $\|x_{k+1}\| = \|x_{k+1}/Y\| = 1$ (this is
possible because $Y$ is finite-dimensional); find
$f_{k+1}' \in (X/Y)^*$ such that $\|f_{k+1}'\| = f_{k+1}'(x_{k+1}/Y) = 1$;
and finally let $f_{k+1}$ be the composition of $f_{k+1}'$ with the natural
projection of $X$ onto $X/Y$. The induction can then proceed.
\medskip
 
Next, make $\N$ into a graph by putting an edge between $k$ and $n$,
$k < n$, if $f_k(x_n) \geq 0$. By infinite Ramsey theory [4, Lemma 29.1]
there is an infinite subgraph which is either complete or anti-complete.
In other words, by passing to a subsequence we can ensure that either
$f_k(x_n) \geq 0$ for all $k < n$ or $f_k(x_n) < 0$ for all $k < n$.
In the first case define $f' = \sum f_k/2^k$ and in the second case
inductively define $a_1 = -1/2$ and
$$a_{n+1} = -{\rm min}(2^{-(n+1)},
{1\over 2} \sum_{k = 1}^n a_k f_k(x_{n+1})),$$
and set $f' = \sum a_k f_k$. It can be seen inductively that each $a_k$ is
negative, and therefore $\sum_1^{n-1} a_k f_k(x_n) > 0$ automatically, and
$a_n$ is chosen to be small enough that also $\sum_1^n a_k f_k(x_n) > 0$.
Thus, in either case, for any $n > 1$ we have $f'(x_n) > 0$.
\medskip

Let $f$ be a weak* cluster point of the sequence $(f_n)$. Then some subnet
of $(f_n)$ converges to $f$, and $f_n$, $f$, and $f'$ all belong to the
unit ball of $X^*$. However, for any $n$ we have
$$(f_n + (f' - f)/2)(x_n) = 1 + f'(x_n)/2 > 1,$$
and since $\|x_n\| = 1$ this implies that $\|f_n + (f' - f)/2\| > 1$, i.e.\
$f_n + (f' - f)/2$ does not belong to the unit ball. Thus, the unit ball
of $X^*$ is not LNC.\hfill\hal
\bigskip

\noindent {\bf Lemma 22.} {\it Every compact, convex, symmetric subset
of a TVS $X$ is linearly homeomorphic to the unit ball of some dual Banach
space, equipped with the weak* topology.}
\medskip

\noindent {\it Proof.} Let $K \subset X$ be such a set and without
loss of generality suppose $K$ spans $X$. Let $V$ denote the space of
linear functionals on $X$ whose restriction to $K$ is continuous, and give
$V$ the sup norm it inherits as a subspace of $C(K)$. Since $V$ is a closed
subspace of $C(K)$, it is a Banach space.
\medskip

Consider the map $T: K \to [V^*]_1$ from $K$ into the unit ball of the
Banach space dual of $V$, defined by $(Tx)(v) = v(x)$. It is easy to check
that this map is continuous going into the weak* topology on $V^*$
and that it is the restriction of a linear map from $X$ to $V^*$.
It is also 1-1 since local convexity of $X$ implies that the continuous linear
functionals on $X$ separate points [5, Corollary 1.2.11]. Thus $K$
is linearly homeomorphic to a weak* compact convex subset of $[V^*]_1$.
\medskip

Now suppose $y \in V^*$ is not in the image of $K$.
Then by a standard separation theorem [5, Theorem 1.2.10] there exists
a weak* continuous linear functional $F$ on $V^*$ such that
$$F(Tx) \leq 1 < F(y)$$
for all $x \in K$. Since $F$ is weak* continuous,
it is again standard [5, Proposition 1.3.5] that there
then exists $v \in V$ such that $z(v) = F(z)$ for all $z \in V^*$, so we have
$$v(x) \leq 1 < y(v)$$
for all $x \in K$. Since $K$ is symmetric and $v$ is linear, this implies
that $\|v\| \leq 1$, hence $\|y\| > 1$.
We conclude that $T$ maps $K$ onto $[V^*]_1$, which completes
the proof.\hfill\hal
\bigskip

\noindent {\bf Theorem 23.} {\it Let $Q$ be a compact, convex, infinite-dimensional
set in a locally convex TVS $X$. Then $Q$ is not LNC.}
\medskip

\noindent {\it Proof.} Let $Q$ be a compact, convex, infinite-dimensional set in
a locally convex TVS $X$, and suppose $Q$ is LNC. Let $Q' = Q - Q$.  Then $Q'$ is
the image of $Q \times Q$ under the map $T: X\times X \to X$ given by $T(x,y) =
x - y$. The map $T$ is continuous and linear, so $Q'$ is compact and convex. Also,
it is easy to see that $Q \times Q$ is LNC, and a trivial variation on Proposition
8 (replacing sequences with nets) therefore implies that $Q'$ is also LNC. So $Q'$
is a compact, convex, infinite-dimensional, symmetric subset of $X$. By Lemma 22
$Q'$ is linearly homeomorphic to the unit ball of some infinite-dimensional dual
Banach space. But this contradicts Lemma 21. We conclude that $Q$ cannot be
LNC.\hfill\hal
\bigskip

\noindent {\bf Corollary 24.} {\it Let $Q$ be a compact, convex, infinite-dimensional
set in a locally convex TVS $E$. Then there is a nonzero vector $v \in E$ such that
the quotient map $T: E \to E/[v]$ restricts to a non-open map $T|_Q: Q \to T(Q)$ and
the section $g_v$ defined as in Section 1 is not continuous.}
\medskip

\noindent {\it Proof.} Identical to the proofs of Theorem 1 (ii) $\Rightarrow$ (i)
and the reverse direction of Proposition 2.\hfill\hal
\bigskip

\noindent {\bf B. Images of closed LNC sets.} The next result shows that under
familiar hypotheses a continuous linear image of a closed LNC set is LNC. In
general such an image need not be closed, so in most applications this will have
to be checked separately.
\bigskip

\noindent {\bf Theorem 25.} {\it Let $Q$ be a closed LNC subset of a Hilbert
space $X$ and let $T: X \to Y$ be a bounded linear map into another Hilbert
space. Suppose $T$ is left Fredholm. Then $T(Q)$ is LNC.}
\medskip

\noindent {\it Proof.} Let $y, y', y_n \in T(Q)$ and suppose $y_n \to y$. Fix
$x, x', x_n \in Q$ such that $T(x) = y$, $T(x') = y'$, and $T(x_n) = y_n$.
Without loss of generality, we may assume that $x$ belongs to the relative
interior of $T^{-1}(y) \cap Q$; translating by $x$, we may further assume
that $x = 0$. As in the proof of Theorem 16, we may also assume that $T = P$
is an orthogonal projection with finite-dimensional kernel.
\medskip

Since $x = 0$ we also have $y = P(x) = 0$, and hence $P(x_n) = y_n \to 0$.
So the hypotheses of Lemma 15 are satisfied, and we conclude that $P'(x_n) \to 0$
and $P'(x_n) \in Q$ eventually. Thus, since $Q$ is LNC we must have
$P'(x_n) + x'/2 \in Q$ eventually; applying $P$ yields $y_n + y'/2 \in T(Q)$
eventually. This verifies that $T(Q)$ is LNC.\hfill\hal
\bigskip

Even in finite dimensions, one can show that if $Q$ is not closed then $T(Q)$
need not be LNC even if $Q$ is LNC. However, any example of such a $Q$ must
be at least four-dimensional, so the verification is bound to be tedious.
Therefore we simply record an example without proving that it has the desired
properties.
\bigskip

\noindent {\bf Example 26.} Define
$$C = \{(x,y,z) \in \R^3: \sqrt{x^2 + y^2} < z < 1\} \cup \{(0,0,0)\}$$
and
$$Q = \{(\lambda x, \lambda y, \lambda z, 1 - \lambda) \in \R^4:
(x,y,z) \in C\hbox{ and }\lambda \in [0,1)\}.$$
Then $Q$ is LNC but its image under the map $T: \R^4 \to \R^3$ given by
$T(x,y,z,w) = (x+w, y, z+w)$ is not.
\bigskip

If $Q$ is convex and open, however, then it is LNC and so is its
image under any bounded linear map $T$ with closed range. For $T$ is an
open map by the open mapping theorem, and so $T(Q)$ is also convex and open,
and hence LNC.
\bigskip

\noindent {\bf C. Unit balls of $L^p$ spaces.} For $1 < p < \infty$, the unit
ball of any $L^p$ space is strictly convex (with respect to the norm topology)
and hence LNC. However, the positive part of the unit ball is not:
\bigskip

\noindent {\bf Proposition 27.} {\it The positive part of the unit ball
of any infinite-dimensional $L^p$ space ($1 \leq p \leq \infty$) is not LNC.}
\medskip

\noindent {\it Proof.} Let $\{A_n: n = 1, 2, \ldots\}$ be a sequence of disjoint positive-measure
sets. Define $x'$ to be the origin of the $L^p$ space and define $x$ by $x|_{A_n} =
\mu(A_n)^{-1/p}2^{-n}$ and $x = 0$ off of $\bigcup A_n$. For each $k \in \N$ define
$$x_k|_{A_n} = \cases{\mu(A_n)^{-1/p}2^{-n}& if $n \neq k$\cr
0& if $n = k$\cr}$$
and set $x_k = 0$ off of $\bigcup A_n$.
Then $x_k \to x$, but $(x_k + (x'-x)/2)|_{A_k} = -\mu(A_n)^{-1/p}2^{-k-1} < 0$.
Thus $x_k + (x'-x)/2$ does not lie in the positive part of the
unit ball.\hfill\hal
\bigskip

Likewise for the unit ball of any infinite-dimensional $L^1$ space:
\bigskip

\noindent {\bf Proposition 28.} {\it The unit ball of any infinite-dimensional
$L^1$ space is not LNC.}
\medskip

\noindent {\it Proof.} Let $\{A_n: n = 1, 2, \ldots\}$ be a sequence of disjoint
positive-measure sets. Define $x|_{A_n} = \mu(A_n)^{-1}2^{-n}$ and
$$x'|_{A_n} = \cases{\mu(A_1)^{-1}& if $n = 1$\cr
0& if $n > 1$.\cr}$$
Also define
$$x_k|_{A_n} = \cases{\mu(A_n)^{-1}2^{-n}& if $1 \leq n < k$\cr
\mu(A_k)^{-1}2^{-k+1}& if $n = k$\cr
0& if $n > k$.\cr}$$
A short computation shows that $x_k \to x$ but
$\|x_k + (x' - x)/2\|_1 = 1 + 2^{-k} > 1$.\hfill\hal
\bigskip

The corresponding fact for $L^\infty$ spaces will be given in Corollary 30.
\bigskip

\noindent {\bf D. Unit balls of $C(K)$ spaces.} For these spaces, unit balls
and their positive parts can be treated simultaneously.
\bigskip

\noindent {\bf Proposition 29.} {\it Let $K$ be a compact Hausdorff space
and suppose $C(K)$ is infinite-dimensional. Then neither the unit ball nor
the positive part of the unit ball of $C(K)$ is LNC. For $K$ locally compact
and $C_0(K)$ infinite-dimensional, the positive part of the unit ball of
$C_0(K)$ is not LNC.}
\medskip

\noindent {\it Proof.} First consider the case that $K$ is locally
compact. Let $(U_n)$ be a sequence of disjoint, nonvoid open
subsets of $K$. For each $n$ let $f_n$ be a Urysohn function supported in
$U_n$. Define a function $g = \sum f_n/2^n$. Define $g_k = \sum_{n \neq k}
f_n/2^n$ and $g' = 0$. Then for $x$ in $U_k$ such that $f_k(x) = 1$,
$$g_k(x) + (g'(x)-g(x))/2 = - 2^{-k-1} < 0.$$
This shows that the positive part of the unit ball of $C_0(K)$ is not LNC.
\medskip

Specializing to the compact case, we find that the positive part of the
unit ball of $C(K)$ is not LNC; taking $h = 1-g$, $h_k = 1-g_k$, and
$h' = 1-g'$ verifies that the entire unit ball is not LNC.\hfill\hal
\bigskip

\noindent {\bf Corollary 30.} {\it The unit ball of any infinite-dimensional
$L^\infty$ space is not LNC.}\hfill\hal
\bigskip

\noindent {\bf Corollary 31.} {\it The positive part of the unit ball
of any infinite-dimensional C*-algebra is not LNC; if the C*-algebra is
unital, its unit ball is also not LNC.}
\medskip

\noindent {\it Proof.} By [1, p.\ 314] every infinite-dimensional
C*-algebra $\A$ contains an infinite-dimensional abelian subalgebra, i.e.\ a
copy of some $C_0(K)$. The positive part of the unit ball of $C_0(K)$ is the
intersection of the positive part of the unit ball of $\A$ with the subalgebra. Since the
intersection of two LNC sets is LNC, and the subalgebra is clearly LNC because
it is a subspace, it
follows that the positive part of the unit ball of $\A$ cannot be LNC.
Likewise for the unit ball of $\A$ if $\A$ has a unit.\hfill\hal
\bigskip

The unit ball of the sequence space $c_0$ is LNC [9, Theorem 1.1], so
Proposition 29 is sharp.
\bigskip
\bigskip

\noindent {\bf 6. Zonoids and related sets.}
\bigskip

A {\it zonoid} is the range of a nonatomic vector-valued measure. By a
well-known theorem of Lyapunov (see e.g.\ [2]), every zonoid is compact
and convex. Equivalently, zonoids are those sets which arise as images of the
positive part of the unit ball of $L^\infty[0,1]$ under weak*-continuous
linear maps into $\R^n$.
\medskip

In this section we will prove that every zonoid is LNC. This is our most
sophisticated construction of finite-dimensional sets with the LNC property,
and it actually can be applied to a somewhat broader class of sets than
zonoids. The basic theorem is the following.
\bigskip

\noindent {\bf Theorem 32.} {\it Let $Q$ be a compact, convex subset of a real
TVS $X$. Suppose that for any closed face $F$ of $Q$ there are compact, convex
sets $A, B \subset X$ such that $Q = A + B$ and $A$ is a translate of $F$.
Then the image of $Q$ under any continuous map $T: X \to \R^n$ is LNC.}
\medskip

\noindent {\it Proof.} Fix $z, z', z_n \in T(Q)$ such that $z_n \to z$ and
let $F'$ be the smallest face of $T(Q)$ which contains $z$ and $z'$.
Now $F = T^{-1}(F') \cap Q$ is a closed face of $Q$, so by hypothesis
there exist compact, convex sets $A, B \subset X$ and $w \in X$ such that
$Q = A + B$ and $F = A + w$.
\medskip

We claim that $a \in A$, $b \in B$,
$a + b \in F$ implies $b = w$. To prove this let $f: \R^n \to \R$ be any
linear functional and $\alpha$ any real number such that
$F' \subset f^{-1}(\alpha)$. Then
$F \subset (f\circ T)^{-1}(\alpha)$ and
$$f(T(a)) = f(T(a + w)) - f(T(w)) = \alpha - \beta$$
for any $a \in A$, where $\beta = f(T(w))$. Therefore, if $a \in A$,
$b \in B$, and $a + b \in F$ then
$$f(T(b)) = f(T(a + b)) - f(T(a)) = \beta,$$
and hence $f(T(A + b)) = \alpha$, so that $T(A + b) \subset
f^{-1}(\alpha)$. Now since $F'$ is a
face of $T(Q)$, for any $u \in T(Q) - F'$ there exists a linear functional
$f: \R^n \to \R$ and a real number $\alpha$ such that
$F' \subset f^{-1}(\alpha)$ and $u \not\in f^{-1}(\alpha)$.
Therefore, by the preceding, $u \not\in T(A + b)$, and we conclude that
$T(A + b) \subset F'$. That is, $A + b \subset F = A + w$. But a compact set
cannot contain a nonzero translate of itself, so $b = w$ as claimed.
\medskip

For each $n$ write $z_n = x_n + y_n$ where $x_n \in T(A)$
and $y_n \in T(B)$, and pass to a subsequence so that $(x_n)$ and $(y_n)$
converge, $x_n \to x$ and $y_n \to y$. Thus $z = x + y$. By the previous
paragraph, this implies that $y = T(w)$. Hence $T(A) = F' - y$.
\medskip

Set $v = (z' - z)/2$.
Since $F'$ is the smallest face which contains $z$ and $z'$, the point
$z + v$ must belong to the interior of $F'$. Subtracting $y$, this
shows that $x + v$ is in the interior of $T(A)$. Therefore $x_n + v \in T(A)$
for sufficiently large $n$, since the convexity notion of interior
coincides with the topological notion inside the span of $T(A)$. Thus
$$z_n + v = x_n + y_n + v \in T(A) + T(B) = T(Q)$$
for sufficiently large $n$, as desired.\hfill\hal
\bigskip

\noindent {\bf Corollary 33.} {\it Let $K$ be a compact, strictly convex
subset of $\R^n$, let $\Omega$ be a $\sigma$-finite measure space, and let
$Q = L^\infty(\Omega; K) \subset L^\infty(\Omega; \R^n)$. Then any
finite-dimensional weak* continuous linear image of $Q$ is LNC.}
\medskip

\noindent {\it Proof.} By Theorem 32 it will suffice to show that $Q$ has
the decomposability property described there. We claim that for any weak*
closed face $F$
of $Q$ there is a subset $\Omega' \subset \Omega$ and a measurable function
$f$ from $\Omega'$ into the boundary of $K$, such that $F$ consists of
precisely those functions in $Q$ which agree with $f$ almost everywhere
on $\Omega'$. From this claim it immediately follows that $Q = A + B$ where
$$A = L^\infty(\Omega - \Omega'; K)$$
is a translate of $F$ (namely, $F = A + f$) and
$$B = L^\infty(\Omega'; K),$$
which verifies the needed decomposition property.
\medskip

To prove the structure theorem, for each $g \in F$ let $\Omega_g$ be the
set of points on which $g$ takes values in the boundary of $K$. Then
let $\Omega'$ be the set whose characteristic function $\chi_{\Omega'}$ is
the infimum in $L^\infty(\Omega)$ of the characteristic functions
$\chi_{\Omega_g}$. It follows that all $g \in F$ take values in the
boundary of $K$ almost everywhere on $\Omega'$, and any $g_1, g_2 \in F$ agree
almost everywhere on $\Omega'$ since otherwise $(g_1 + g_2)/2$ would take
values in the interior of $K$ on a positive measure subset of $\Omega'$,
contradicting its definition. Thus we define $f = g|_{\Omega'}$ where
$g$ is any function in $F$.
\medskip

Now let $h$ be any function in $Q$ which agrees with $f$ almost everywhere
on $\Omega'$. We will show that $h$ is a weak* limit of functions in $F$,
hence $h \in F$. This will complete the proof. For any $g \in F$ and any
$\epsilon > 0$ the function $h_\epsilon$ defined by
$$h_\epsilon(x) = \cases{g(x)& if $d(g(x), \R^n - K) \leq \epsilon$\cr
h(x)& otherwise\cr}$$
belongs to $F$ since
$$g = \lambda h_\epsilon +
(1 - \lambda)(g + {\lambda\over{1 - \lambda}}(g - h_\epsilon))$$
and $g + {\lambda\over{1 - \lambda}}(g - h_\epsilon) \in F$ for sufficiently
small $\lambda$. Taking the limit of $h_\epsilon$
as $\epsilon \to 0$, we find that $F$
contains the function $h_g$ which agrees with $g$ when $g(x)$ belongs to
the boundary of $K$ and agrees with $h$ when $g(x)$ belongs to the interior
of $K$. Now for any finite set $\{g_1, \ldots, g_n\} \subset F$ the average
$g = (g_1 + \cdots + g_n)/n$ also belongs to $F$, and $g(x)$ belongs to
the boundary of $K$ only if $g_i(x)$ belongs to the boundary of $K$ for all
$i$. Thus $h$ is a weak* cluster point of the functions $h_g$ as $g$ ranges
over $F$, and we conclude that $h \in F$.\hfill\hal
\bigskip

\noindent {\bf Corollary 34.} {\it Every zonoid is LNC.}
\medskip

\noindent {\it Proof.} Take $\Omega = [0,1]$ and $K = [0,1] \subset \R$ in
Corollary 33. Then $Q$ is the positive part of the unit ball of
$L^\infty[0,1]$, and by Corollary 33 any weak* continuous linear image of
$Q$ is LNC.\hfill\hal
\bigskip

\noindent {\bf Example 35.}
Note that a continuous image of the positive part of the unit ball of $l^1$
need not be LNC. For example, let $C$ be the convex hull of the origin in
$\R^3$ and the points $p_n = (1/n, 1/n^2, 1)$ ($n \in \N$). This set is not
LNC by Proposition 11. Define $T: l^1 \to \R^3$ by $T(f) = \sum f(n)p_n$.
Then the image under $T$ of the positive part of the unit ball of $l^1$
equals $C$, which is not LNC.
\bigskip

\noindent {\bf Example 36.}
Similarly, the positive part of the unit ball of any nonabelian von Neumann
algebra has weak* continuous, finite-dimensional, linear images which are
not LNC. To see this, let $\M$ be a nonabelian von Neumann algebra and
find a subalgebra $\A$ which is isomorphic to $M_2(\C)$ [10, p.\ 302].
Find a faithful normal representation of $\M$ on a Hilbert space $H$.
Since any representation of $M_2(\C)$ can be decomposed into an orthogonal
direct sum of two-dimensional representations, we can find a two-dimensional
subspace $K$ of $H$ such that $p\A p = B(K)$ where $p$ is the orthogonal
projection of $H$ onto $K$. Now define $T: \M \to B(K)$ by $Tx = pxp$. This is
a weak*-continuous linear map, and it takes the positive part of the unit
ball of $\M$ onto the positive part of the unit ball of $B(K) \approx
M_2(\C)$, which is not LNC by Proposition 12.
\bigskip

\noindent {\bf Example 37.}
We also note that there is a finite-dimensional subspace of $l^\infty$ whose
intersection with the positive part of the unit ball is not LNC. To prove
this let $D$ be a countable set of unit vectors in $\C^2$, dense in the set
of all unit vectors, and for each
$x_n \in D$ ($n \in \N$) define $f_n: M_2(\C) \to \C$ by $f_n(A) =
\langle Ax_n, x_n\rangle$. Now define $T: M_2(\C) \to l^\infty$ by $T(A) =
(f_n(A))$ (the sequence whose $n$th term is $f_n(A)$).
This map is nonexpansive because each $f_n$ has norm one, and in fact
$\|T(A)\| = \|A\|$ for all positive matrices $A$.
Thus $T(M_2(\C))$ is a finite-dimensional subspace of
$l^\infty$ whose intersection with the positive part of the unit ball of
$l^\infty$ is isometric to the positive part of the unit ball of
$M_2(\C)$, which again is not LNC by Proposition 12.
\bigskip
\bigskip

\noindent [1] C.\ A.\ Akemann, Left ideal structure of C*-algebras,
{\it J.\ Funct.\ Anal.\ \bf 6} (1970), 305-317.
\medskip

\noindent [2] C.\ A.\ Akemann and J.\ Anderson, {\it Lyapunov Theorems for
Operator Algebras}, {\it Mem.\ Amer.\ Math.\ Soc.\ \bf 458} (1991).
\medskip

\noindent [3] C.\ A.\ Akemann and N.\ Weaver, Continuous selection of
bang-bang controls, unpublished manuscript.
\medskip

\noindent [4] T.\ Jech, {\it Set Theory}, Academic Press (1978).
\medskip

\noindent [5] R.\ V.\ Kadison and J.\ R.\ Ringrose, {\it Fundamentals of the Theory     
of Operator Algebras II}, Academic Press (1986).
\medskip

\noindent [6] E.\ Michael, Continuous selections I,
{\it Ann.\ of Math.\ \bf 63} (1956), 361--382.
\medskip

\noindent [7] W.\ Rudin, {\it Functional Analysis} (second edition),
McGraw-Hill (1991).
\medskip

\noindent [8] G.\ C.\ Shell, {\it Locally Nonconical Convex Sets},
Ph.D.\ Thesis, UC Santa Barbara (1998).
\medskip

\noindent [9] ---------, On the geometry of locally nonconical convex sets,
{\it Geom.\ Dedicata \bf 75} (1999), 187-198.
\medskip

\noindent [10] M.\ Takesaki, {\it Theory of Operator Algebras I},
Springer-Verlag (1979).
\bigskip
\bigskip

\noindent Department of Mathematics

\noindent University of California

\noindent Santa Barbara, CA 93106

\noindent akemann@math.ucsb.edu
\bigskip

\noindent Department of Defense

\noindent 9800 Savage Road

\noindent Fort Meade, MD   20755-6000

\noindent GCShell@aol.com
\bigskip

\noindent Department of Mathematics

\noindent Washington University

\noindent St.\ Louis, MO 63130

\noindent nweaver@math.wustl.edu
\end